\input amstex

\documentstyle{amsppt} \magnification1200
\TagsOnRight
\pagewidth{30pc} \pageheight{47pc}

\input epsf.tex
\def\today{1 Sept 2000}

\def\R{{\Bbb R}}

 \def\epsfcenter#1{{\vcenter{\hbox{\epsfbox{#1}}}}} 
\def\epsfpad#1#2 #3 {{\vcenter{\vskip#2\hbox{\epsfbox{#1}}\vskip#3}}} 
\topmatter
\title Unlinked Embedded Graphs\endtitle
\author John W. Barrett\endauthor
\date\today\enddate
\address
School of Mathematical Sciences,
University of Nottingham,
University Park,
Nottingham,
NG7 2RD, UK
\endaddress
\email jwb\@maths.nott.ac.uk \endemail
 
\abstract This paper is a self-contained development of an invariant of graphs embedded in three-dimensional Euclidean space using the Jones polynomial and skein theory. Some examples of the invariant are computed. An unlinked embedded graph is one that contains only trivial knots or links. Examples show that the invariant is sufficiently powerful to distinguish some different unlinked embeddings of the same graph. 
\endabstract              

\endtopmatter
 
\document


The Jones polynomial\cite{J} assigns an invariant to each oriented link $L$ in $\R^3$ which is a Laurent polynomial $V_L(t)$ in the variable $\sqrt t$. The product $V_L(t)V_L(t^{-1})$ is another Laurent polynomial in $t$ and can be expressed in terms of $\lambda=t+t^{-1}$ to give an ordinary polynomial

$$r_L(\lambda)=r_L(t+t^{-1})= V_L(t)V_L(t^{-1}).$$
The polynomial $R_L(\lambda)$ is defined by
\footnote{The extra factor $\lambda+2$ is included so that the unknot $U$ has invariant $R_U(\lambda)=\lambda+2$ and the empty link has polynomial $1$.}
$$R_L(\lambda)=(\lambda+2)r_L(\lambda).$$

The invariant $R_L$ does not depend on the orientation of the link $L$ and cannot distinguish a link from its mirror image. The advantage of the polynomial $R_L$ is that the definition extends to an invariant of graphs embedded in $\R^3$, a generalisation of the idea of a link in which vertices where a number of edges meet are allowed. The idea of this paper is to develop this invariant using elementary ideas from skein theory and this definition of $R_L$ for links. The paper is restricted to the simplest case for which every vertex is four-valent (has degree four), but there is a generalisation to vertices of any even valence.

A graph is defined here as a compact polyhedron which is locally isomorphic to either an interval of the real line or to $\epsfcenter{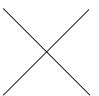}$; in particular a component of a graph may be a circle with no vertices.

The invariant for the four-valent vertex is defined by
$$\epsfcenter{vertex.eps}
= P(\lambda)\left( \quad
\lambda\;\epsfcenter{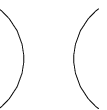}+\lambda\;\epsfcenter{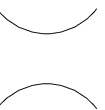}
+\epsfcenter{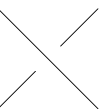}+\epsfcenter{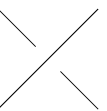}
\quad \right)
\tag1$$
The graph invariant $R_G$ for an embedded graph $G$ is defined by applying this relation to every vertex in a diagram for $G$, then evaluating the resulting link diagrams $L$ using the previously defined $R_L$. The fact that this does not depend on the projection and gives an invariant of ambient isotopy of the embedding is proved below.

The normalising factor $P(\lambda)$ can be chosen arbitrarily. Choosing $P=1$ would give a polynomial for each embedded graph. The definition used here is
$$ P(\lambda)={1\over(\lambda+1)(\lambda+2)}$$
This choice has the disadvantage that the invariant for graphs is no longer always a polynomial but a rational function. However the advantage is that the expressions for the examples calculated here are somewhat simpler.

Louis Crane and I came across a formula for the 4-valent vertex while studying quantum gravity, expressed in terms of the representation theory of $U_qsl(2)$ \cite{BC}. Yetter gave the theory of the embedded graph invariant, and extended it to vertices of arbitrary valence\cite{Y2}. In \cite{B} I gave an alternative definition of Yetter's invariant based on the Kauffman bracket approach to spin networks \cite{K2}. The invariant in those papers is more general than the one considered here in that each edge is labelled with an arbitrary representation of $U_qsl(2)$. The idea of the present paper is to take a special case, namely the fundamental representation, and to give a completely self-contained treatment which does not need any machinery of quantum groups or category theory. 
In the process, simple formulae, such as the definition \thetag1, become apparent, and the properties and examples can be developed rapidly.

\topinsert \topcaption{Table 1} The embedded graph invariant\endcaption
$$\matrix
&&R(\lambda) & R(2)& R(1)& R(-1)\\
A&\epsfpad{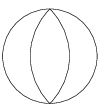}5pt 5pt &2&2&2&2\\
A'&\epsfpad{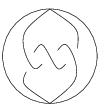}5pt 5pt &-\lambda^3+2\lambda^2+2&2&3&5\\
B&\epsfpad{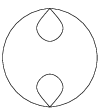}5pt 5pt &\lambda+2&4&3&1\\
B'&\epsfpad{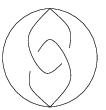}5pt 5pt &\lambda^2-\lambda+2&4&2&4\\
B''&\epsfpad{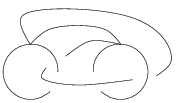}5pt 5pt &-\lambda^3+4\lambda^2-3\lambda+2&4&2&10\\
C&\epsfpad{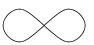}5pt 5pt &\lambda+2&4&3&1\\
C'&\epsfpad{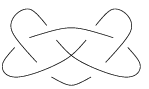}5pt 5pt & (\lambda+2)(\lambda^5-3\lambda^4+4\lambda^2+1) &4 &9 &1 \\
D&\epsfpad{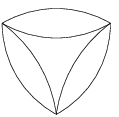}5pt 5pt & {\lambda+2\over\lambda+1} &{4\over3} &{3\over2} & \infty \\
D'&\epsfpad{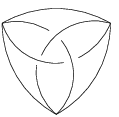}5pt 5pt & {-\lambda^3+3\lambda^2-\lambda+2\over\lambda+1}& {4\over3} &{3\over2} &\infty \\
E&\epsfpad{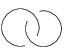}5pt 5pt & \lambda^2(\lambda+2)& 16&3 &1 \\
F&\epsfpad{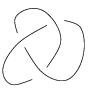}5pt 5pt & (\lambda+2)(-\lambda^3+\lambda^2+2\lambda+1) & 4&9 & 1\\
G&\epsfpad{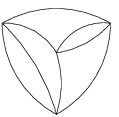}5pt 5pt &{3\lambda^2+6\lambda+4\over(\lambda+1)^2(\lambda+2)} &{7\over9} &{13\over12} & \infty \\
\endmatrix$$
\endinsert

Table 1 gives some examples of the evaluation of the invariant for 4-valent graphs and links.
In the table, the examples with the same letter, e.g. $A$ and $A'$, are the same graph but with different embeddings. Two edges can be removed from $A'$ to give a trefoil knot $F$ and so it is perhaps not surprising that the invariant can distinguish $A$ and $A'$. The examples $B$,$B'$,$B''$ again share the same graph. $B'$ is {\it linked} in the sense that removing the two outer edges gives the Hopf link $E$. 

However  $B''$ is an {\it unlinked embedded graph} in this sense: any way of removing edges from $B''$ to make a link results only in a number of unlinked unknots. Yet $B$ and $B''$ differ; they have different $R$ invariants.
The example $C'$ is also an unlinked embedding\footnote{This example was found by Paul Langlois} which differs from the `trivial' embedding $C$.

\subhead Properties of the invariant \endsubhead
The main relation satisfied by the invariant $R_L$ for links is the following cubic relation on the braid generator

\proclaim {Lemma 1} The invariant satisfies
$$\epsfcenter{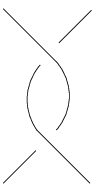}= 
(1-\lambda) \epsfcenter{crossing.eps}
+(\lambda-1)\epsfcenter{vertical.eps}
+\epsfcenter{othercrossing.eps}$$\endproclaim

\demo{Proof} 
The invariant $R_L$ for a link $L$ is related to the (suitably normalised) Kauffman bracket polynomial\cite{K} for the same diagram $N(A)$, by
$$ R_L(A^4+A^{-4})=N(A)N(A^{-1}).$$
The braid generator $b$ for the Kauffman bracket satisfies the quadratic relation
$$A^{-1} b^2 + (A^2-A^{-2})b-A=0.$$
The braid generator for the invariant $R_L$ can be represented by
$$\epsfcenter{crossing.eps}=b\otimes b^{-1}.$$ Then Lemma 1 follows from using the quadratic relation for $b$.\enddemo

\proclaim{Theorem 2}
The definition of $R$ for 4-valent embedded graphs in $\R^3$ is independent of the diagram and is an invariant of ambient isotopy.\endproclaim

To prove invariance under ambient isotopy of the embedding of the graph it suffices to check a set of extended Reidemeister moves \cite{Y}.
The symmetry relation
$$\epsfcenter{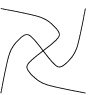}=\epsfcenter{vertex.eps}$$
follows immediately from the definition, while the permutation property of the vertex
$$\epsfcenter{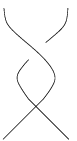}=\epsfcenter{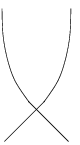} 
\tag2$$
follows from a calculation using Lemma 1.
\bigbreak
This theorem could also be proved by establishing the equivalence with the invariants in \cite{Y2}.

The invariant is an example of the more general `rigid vertex' invariants described in \cite{KV} which do not necessarily have the invariance \thetag2 under the permutations of edges meeting at the vertex. Some other equivalence relations on embedded graphs and invariants of these are studied in \cite{T}. 

\bigbreak

\proclaim{Lemma 3}
A further property of the invariant is
$$ \epsfcenter{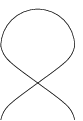}=\epsfcenter{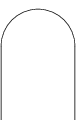}
$$
\endproclaim
This can be taken as the way of fixing the normalising factor $P$.

\subhead  Specializations\endsubhead
Finally, there are three specializations of particular interest, $\lambda=2$, $1$ and $-1$, tabulated for the examples in table 1. The value  $R(2)$ does not depend on the embedding of the graph. The invariant is just a product of the invariants for each component, with a circle equal to 4. The invariant $R(1)$ for links gives the three-colouring invariant \cite{P}. For a graph, this gives an extension of the three-colouring invariant as a count of the number of ways of colouring the arcs of a diagram with particular weights for the patterns at each vertex. The invariant $R(-1)$ is equal to 1 for every link, but this no longer holds for embedded graphs.

\subhead Acknowledgement \endsubhead
 Thanks are due to the hospitality of the Instituto Superior T\'echnico 
 Lisbon in December 1999.

 \enddocument